\documentclass[11pt]{article}
\usepackage{amsmath,amssymb}
\usepackage{enumerate}
\usepackage{color}
\usepackage{subfigure}
\usepackage{tikz-cd}

\newcommand{\R}{\mathbb{R}}

\newcommand{\N}{\mathbb{N}}
\newcommand{\bea}{\begin{eqnarray}}
	\newcommand{\eea}{\end{eqnarray}}

\def\a{\alpha}

\def\e{\varepsilon}

\def\s{\sigma}

\def\x{\otimes}

\def\supp{{\rm supp}}

\def\ric{{\rm ric}}
\def\1{\rm Id}

\def\ci{\circ}

\def\pr{{\rm pr}}
\newcommand{\qed}{$\hfill\blacksquare$}

\def\V{\noindent}

\def\grad{{\rm grad}}

\def\Int{{\rm int}}

\newcommand{\bean}{\begin{eqnarray*}}
	\newcommand{\eean}{\end{eqnarray*}}

\newtheorem{Theorem}{Theorem}

\newcommand{\ben}{\begin{enumerate}}
	\newcommand{\een}{\end{enumerate}}
\newcommand{\bit}{\begin{itemize}}
	\newcommand{\eit}{\end{itemize}}
\newcommand{\edoc}{\end{document}}

\usepackage{hyperref}

\usepackage{enumerate}

\title{Bartnik's splitting conjecture with the null energy condition}

\begin{document}

\author{Olaf M\"uller\footnote{Institut f\"ur Mathematik, Humboldt-Universit\"at zu Berlin, Unter den Linden 6, D-10099 Berlin, \texttt{Email: o.mueller@hu-berlin.de}}}

\date{\today}
\maketitle

\begin{abstract}
	\V Bartnik's splitting conjecture is one of the prime open conjectures in mathematical relativity. There are many approaches to this conjecture that use (Lorentzian) conformal geometry. In this article, we show that if we replace the strong energy condition in Bartnik's splitting conjecture with the null energy condition, then in any dimension greater or equal to $3$ the conclusion of the conjecture would be wrong, more precisely: On a manifold of dimension at least $ 3$, {\em every} globally hyperbolic spatially compact conformal class contains future complete metrics satisfying the null energy condition. In the spatially noncompact case, the same is true in the future of any Cauchy surface. The main tool is the flatzoomer method.
\end{abstract}

\section{Introduction and result}

\V Bartnik's famous splitting conjecture \cite{rB1} --- inspired by the idea of exploring rigidity of the Penrose-Hawking singularity results --- is the following:

\bigskip

\V {\bf Conjecture:} Any $n$-dimensional timelike geodesically complete, globally hyperbolic spatially compact\footnote{i.e., with a compact Cauchy surface} manifold $(M,g)$ that is {\em SEC} (i.e.,  satisfies the strong energy condition $\ric_g (v,v) \geq 0$ for each $g$-causal $v \in TM$) {\em splits}, i.e., is isometric to a semi-Riemannian product of an $(n-1)$-dimensional  Riemannian manifold and the Lorentzian line. 

\bigskip

\V Several partial results on this topic have been obtained during the last decades due to Bartnik, Galloway, Gerhardt and others. For example, the conjecture is true for $\dim (M)=2$, as in this case the energy condition is equivalent to sectional curvature being nonnegative. Bartnik's conjecture is also true under Galloway's ray-to-ray condition \cite{Greg-ray-to-ray}, and also if we assume the existence of a {\em CMC} Cauchy hypersurface (i.e., one of constant mean curvature), cf. Cor. 2 in \cite{rB1}, and Bartnik then in \cite{rB1} showed this to be implied by {\em Bartnik's conformal condition (BCC)}: that there be some $p \in M$ with $M \setminus J(p)$ precompact, where $J(p) := J^+(p) \cup J^-(p)$ --- which is equivalent to saying that $J^+(p) $ and $J^-(p)$ both contain a Cauchy surface\footnote{Given a Cauchy temporal function $t$ on $X$, the first implication follows from $t$ taking a minimum $T_-$ and a maximum $T_+$ on $ X \setminus J(p)$, then $ t^{-1}(\{ T_\pm \pm 1 \})$ satisfy the condition. The other implication follows from compactness of $ J(U,S)$ for Cauchy surfaces $U,S$ in a spatially compact spacetime.}. Spacetimes satisfying this BCC condition will henceforth be called {\bf long}. This is a conformally invariant condition, singling out a subset of the set of conformal structures. 

It is probably safe to say that all {\em purely causal} conditions on a cosmological SEC and timelike complete spacetime under which the conjecture is known to be true are special cases of the {\em min-min condition} presented by Leonardo Garcia-Heveling \cite{lGH}: that $M$ contains a TIP and a TIF (terminal indecomposable past resp. future sets) both minimal w.r.t. the order in the TIPs/TIFs by inclusion and intersecting nontrivially. Here, a terminally indecomposable past set is a past set $A$ in $M$ (i.e., $I^-(A) = A$) that is not the past of a point (i.e., $A \neq I^-(p) \ \forall p \in M$) and is only trivially decomposable into past sets, i.e., if $A = B \cup C$ for past sets $B,C$ then $B \subset C $ or $C \subset B$. A classical theorem (\cite{GKP}, Th.2.1) states that if a spacetime $X$ is distinguishing (which is implied by our general assumption of $X$ being globally hyperbolic) for each TIP $A$ in $X$ there is a $C^0$ inextendible future curve $c: [0; \infty) =:K \rightarrow M$ such that $A = I^-(c(K))$, and a corresponding time-dual statement holds for TIFs.

\bigskip

\V This article considers a weakening of the hypotheses from SEC to NEC (null energy condition). De Sitter space shows already that the NEC is not sufficient for Bartnik's conjecture. Here we show that, even more, {\em in dimension $\geq 3$, the combination of future\footnote{Compare with the recent result by Galloway and Ling \cite{GL2} that each future timelike complete cosmological spacetime satisfying the timelike convergence condition with a smooth Cauchy surface of nonpositive mean curvature contains a maximal Cauchy surface.} timelike (and causal) completeness and the NEC does not restrict the conformal structure of the spacetime at all}. 

\bigskip

For a real interval $I$, an arclength-parametrized timelike curve $c \in C^2(I, M)$ in a spacetime $(M,g)$ is called {\bf of bounded acceleration} or {\bf b.a.} iff $t \mapsto g(\nabla_t (c'(t)), \nabla_t (c'(t))) $ is bounded on $I$. A spacetime is called {\bf future (resp. past) b.a. complete} iff each future (resp. past) $C^0$-inextendible b.a. curve has infinite length (or, equivalently, infinite parameter), cf \cite{BEE} Sec. 6. We get:

\begin{Theorem}
	\label{main-flatzoomer}
	Let $(M,g)$ be a globally hyperbolic Lorentzian manifold of dimension $\geq 3$. Then: 
	\begin{enumerate}
		\item For every spacelike Cauchy surface $N$ of $(M,g)$, there is $u \in C^\infty (M)$ such that $(M, e^{2u} \cdot g)  $ is causally future geodesically and b.a.-complete, and that satisfies the NEC on $J^+(N)$.
		\item If $(M,g) $ is spatially compact, then there is $u \in C^\infty (M)$ such that $(M, e^{2u} \cdot g)  $ is causally future geodesically and b.a.-complete and satisfies the NEC.
	\end{enumerate}
\end{Theorem}

\bigskip

{\bf Conventions of the article:} $0 \in \N$, $\N^* := \N \setminus \{0\}$, $\N_n:= \{ k \in \N | k \leq n \}$, $\N_n^*:= \N_n \cap \N^*$. Intervals are written with semicola instead of commata, to distinguish open intervals from tuples.
\newpage

\section{Proof of the main result}
\label{SECNEC}

The result is a combination of two earlier results of the author (\cite{EMBH}, Th. 16, \cite{oM:Horizons}, Th.8), stating existence of one conformal factor satisfying both requirements at once. Its proof, however, requires methods different from those in \cite{EMBH} and \cite{oM:Horizons}. For a subset $A$ of a spacetime $X$ we define $\partial^\pm A:= \{ x \in \partial A | I^\pm(x) \cap A = \emptyset \}= \partial A \cap (X \setminus I^\mp (A))$, a closed subset of $X$.

\V A function $t$ on a spacetime $(M,g)$ is called {\bf Cauchy} iff for any $C^0$-inextendible timelike curve $c: (T_-; T_+) \rightarrow M$, the function $t \circ c: (T_-; T_+) \rightarrow \R$ is surjective, and {\bf temporal} iff it is $C^1$ and $\grad (t)$ is past timelike, and {\bf steep} iff it is temporal and $g(\grad t, \grad t) \leq -1$. Temporalness ensures injectivity of $t \circ c$ for all timelike curves $c$, thus the level sets of a Cauchy temporal function are Cauchy surfaces. Each globally hyperbolic (g.h.) spacetime admits a steep Cauchy temporal function vanishing on a given Cauchy surface (\cite{oM-inv}, \cite{oMmS}). 

\begin{Theorem}
	\label{Lemma2}
	Let $(M,g)$ be globally hyperbolic, let $S \subset (M,g)$ be a spacelike Cauchy hypersurface. Then there is a steep Cauchy function $t$ with $t^{-1} ((- \infty; D)) \cap J^+(S)$ compact $ \forall D\in \R$. 	
\end{Theorem}

\vspace{-0.2cm}

\V {\bf Proof.}  Let $t_0 $ be a smooth steep Cauchy function with $S= t_0^{-1} (\{ 0 \} )$, whose existence is ensured by \cite{oM-inv}[Th.1], let $S_r:= t_0^{-1} (\{ r\})$. For $n \in \N$, let $U_n \subset S$ be compact with $U_n \subset \Int (U_{n+1})$ for all $n \in \N$, and $ \bigcup_{n =0}^{\infty} U_n = S$. By an obvious inductional argument, using compactness of $ J^\pm (K) \cap J^\mp (S)$ for compacta $K$ and Cauchy surfaces $S$, we can assume w.l.o.g. that $\{ U_n | n \in \N\}$ is {\bf safe}, i.e. $V(U_n) := J^-(J^+(J^- (U_n) \cap S_{-1}) \cap S_n) \cap S_0 \subset U_{n+1}$.  The set $K_n:= J^+(U_n) \cap J^-(S_n) $ is compact, and $K_n \subset \Int (K_{n+1}) \cup U_n \subset I^-(K_{n+1})  \ \forall n \in \N$ by continuity of $J^-: M \rightarrow 2^M$ (\cite{BEE} pp.59-73, \cite{oM-FCC}), $ \bigcup_{n \in \N} K_n = J^+(S)$ by Cauchyness of $t_0$ and $S$. Defining $u_k:= U_{k+1} \setminus U_k \forall k \in \N $ we get $J^+(u_k) \cap J^+(u_l) \cap J^-(S_n) \neq \emptyset \Rightarrow |k-l| \leq 1 \forall k,l \geq n $. We define $t: M \rightarrow \R$ steep Cauchy with 

\bea
\label{more}
t\vert_{J^+(S) \setminus  K_{n-1} } \geq n \ \forall n \in \N,
\eea

\V as follows: We have $\partial K_n \cap I^+(S) = (\partial^+ K_n) \cap I^+(S) \subset \Int K_{n+1}$, define $C_n:= K_{n} \setminus \Int (K_{n-1} ) $. For every $p \in C_n$ there are $p_{- -} \ll p_- \ll p $ with $p_{--} \in L_n:= (J^-(U_n \setminus \Int U_{n-1}) \cap J^+(S_{-1})) \cup J^+ (S_{n-2}) $. By compactness of the $C_n$, we find $m \in \N$ and a finite subset $\{ p^{1,n}_{--},..., p^{m,n}_{--} \} \subset L_n$ such that $ \bigcup_{i=1}^m I^+(p^{i,n}_-)$ covers $C_n$. Let $\psi \in C^\infty (\R)$ with $\psi (x) = 0 \ \forall x \leq 0 $, $\psi ' (x) > 0 \ \forall x >0$ and $\psi' (x) \geq 1 \ \forall x \geq 1$. Let $t_{i,n}$ be a smooth steep Cauchy temporal function on $B^{i,n}:= I^+(p^{i,n}_{--})$ (a causally convex subset and thus g.h.) with $ t_{i,n} (p^{i,n}_-) = 1$. We put $\tau_{i,n} := \psi \ci t_{i,n}$ on $B^{i,n}$ and $\tau_{i,n} \vert_{M \setminus B_i^n} = 0$. Then, for $a_{i,n} >0$ sufficiently large, the function $t_n:= \sum a_{i,n} \cdot \tau_{i,n} \in C^{\infty} (M, (0; \infty))$ is, by compactness of $C_n$, a smooth steep temporal function on an open subset containing $C_n$ with $t_n |_{C_n} \geq n$. By safeness of $ \{ U_n | n \in \N\}$ we get $\supp (t_{n+1}) \cap J^+(S)  \subset J^+(U_{n+2} \setminus U_{n-1}) \cup J^+(S_{n-2}) $.
Finally, $t:= t_0 + \sum_{i=1}^{\infty} t_i$ is well-defined and smooth on $M$, as for any fixed point $q$ there are finitely many nonzero contributions from the sum: if $J^-(q) \cap S \subset U_n$ then $q \notin \supp (t_N) \forall N > n+1, t_0(q)+2$. Finally, $t$ is a Cauchy steep function being the sum of the Cauchy steep function $t_0$ and temporal functions $t_i$, see \cite{oMmS}, Prop. 2.1 (c), and it satisfies Eq. \ref{more}. \hfill \qed 

\bigskip

\V {\bf Proof of Theorem \ref{main-flatzoomer}.} We consider first the case that $ (M,g) $ is spatially compact. Let $t$ be a smooth Cauchy temporal function on $M$. By its gradient flow, $t$ induces an isometric diffeomorphism $F: (M,g) \rightarrow (\R \times N , - k \cdot dt^2 + \pr_2^*(g^N (t) ) $ where $\pr_2: \R \times N \rightarrow N $ is the projection to the second component, $k \in C^\infty (\R \times N, (0; \infty))$ and $t \rightarrow g^N (t) $ is a smooth curve of Riemannian metrics on $N$. Using the freedom to choose conformal factors, we can assume w.l.o.g. that $k$ is identically $1$. For any conformal multiple $Ug = e^{2u} g$ of $g$, the very same isometry $F$ induced by the function $t$ is an isometry between $Ug $ and $ - U dt^2 +  U g_t$. The geodesic equation in $(M,Ug)$ for the $t$ coordinate along a geodesic $c$ (where $q' := \frac{\partial}{\partial t} q$ for a differentiable function $q$ on $\R \times N$) reads after calculation of the Christoffel symbol:

\bea
\label{TproS}
\frac{d^2 t}{ds^2}  = - \sum_{m=1}^{n-1}  \frac{U_{,m}}{U} \frac{dx_m}{ds} \frac{dt}{ds} - \frac{U'}{2U}  (\frac{dt}{ds} )^2    - \sum_{m,n =1}^{n-1} \frac{1}{2} (g'_{mn} + \frac{U' }{U} g_{mn} )   \frac{dx_m}{ds}   \frac{dx_n}{ds}, 
\eea

\V where $s$ is the affine parameter of $c$. We want to show that for the appropriate growth of $U$, Eq. \ref{TproS} implies convexity of $t \ci c$. With the ansatz $U = A \ci t$ for some $A = e^{2 a} \in C^\infty (\R, \R)$, the first term vanishes. If we consider a b.a. curve we get a bounded real function $D$ as an additional additive term on the right hand side. For the geodesic case just replace $D$ by $0$ from now on. Because of compactness of the $t$-level sets $S_t$ and positive-definiteness of its Riemannian metric $g^{(t)}$, there is a $\eta(t) >0 $ such that $(g^{(t)})' + 2 \eta(t) \cdot g^{(t)} $ is positive definite. Now, if 

\bea\label{FZ1}
\gamma (t):= \frac{d (\ln (A(t)))}{dt} = 2a'(t)  \geq 2 \eta(t)   ,
\eea

then both remaining terms on the RHS of Eq. \ref{TproS} are nonnegative, thus $\frac{d^2 t}{ds^2} \leq D$, but as $ \frac{ dt}{ds} >0$, by the mean value theorem we get $s \geq (\frac{dt}{ds})^{-1} \vert_{s=0} \cdot t$ (geodesics) or $s \geq D (\frac{dt}{ds})^{-1} \vert_{s=0} \cdot t^2  $ (b.a. curves), so in either case, for bounded $s$, $t$ is bounded as well, thus $(I^+ (S), U \cdot g )$ is future geodesically complete and b.a.-complete.

\V On the other hand, for the energy condition, we consider the well-known equation

\bea
\widetilde{\ric} = \ric - (n-2) (\nabla^g du - du \x_{{\rm sym}} du ) + (\Delta u - (n-2) g(\grad^g u, \grad^g u) ) \cdot g 
\eea

for the Ricci curvature $\widetilde{\ric}$ of the metric $\tilde{g} := e^{2u} g$ (see for example Eq. 1.159 d) from \cite{BESSE}). In our case we recall that $u = (\ln U)/2$ and $U= A \ci t$ for some $A: \R\rightarrow (0; \infty)$. Let $\tau M : TM \rightarrow M$ be the tangent bundle projection, let $K_b$ be the compact set of causal vectors $v$ in $ D_b:= (\tau M)^{-1}  (t^{-1} (b)) $ such that $dt (v) =1$, then every causal vector in $D_b$ is a multiple of a vector in $ K_b$. Focussing on the quadratic terms $du \x_{{\rm sym}} du$ and $g(\grad^g u, \grad^g u) ) \cdot g$ we see that the latter is nonpositive on $K_b \times K_b$ whereas the former is positive on $K_b \times K_b$ if $A '(x) >0 \ \forall x \in \R$, and due to compactness of $K_b$ it is even bounded below by a positive constant. Let $H_b$ be the intersection of $K_b$ with the null cone, then

\bea
\widetilde{\ric} |_{H_b \x H_b}= \ric - (n-2) ( \nabla^g du - du \x_{{\rm sym}} du )|_{H_b \times H_b} ,
\eea

We want to make use of the theory of flatzoomers as defined in \cite{MN}. For $u = a \ci t $ for some $a \in C^1 (\R, \R) $ we get $du = (a' \ci t)  \cdot dt   $, $du \x_s du = (a ' \ci t )^2 \cdot dt^2$,

\[  \nabla^g du = \nabla^g ( (a' \ci t)  \cdot dt)  =  (a'' \ci t) \cdot dt^2 + (a' \ci t) \cdot \nabla^g dt \] 

Now we define $\beta := \ln (a') $, so $a' = e^{\beta} $ and $a'' = \beta ' \cdot e^{\beta} $. We get 

\[\tilde{\ric}|_{H_b \times H_b} = \ric|_{H_b \times H_b} - (n-2) \psi (\beta)|_{H_b \times H_b} , {\rm \ with \ } \psi(\beta) := e^{\beta \ci t} \cdot ((\beta ' \ci t) dt^2 + \nabla^g dt - e^{\beta \ci t} dt^2 ): TM^2 \rightarrow \R .\]

For each $b \in \R$ we define 

\[
\phi(\beta) (b) := \max \{ \psi (\beta) (v,v) \vert v \in H_b \}< \infty, \ \ \rho (b) : = \min \{ \ric (v,v) \vert v \in H_b \}> - \infty,\]

both estimates above by compactness of $H_b$. If we find $\beta$ with $\phi(\beta) (b) < \rho (b)/(n-2) \ \forall b \in \R$, then for $a(s):= \frac{1}{2}\int_0^s e^{\beta (r)} d r$ and $u:= a \ci t$, the metric $\tilde{g} = e^{2u} g$ satisfies $\ric^{\tilde{g}} (v,v) >0 \ \forall v \in H_b \ \forall b \in \R$ (note that $n-2 >0$), so $\tilde{g}$ satisfies the NEC. 

We put $\a (b) := \max \{ \nabla^g dt (v,v) | v \in H_b \} < \infty$. Then we get 

\[ \phi (\beta) (b) = e^{2\beta} (e^{- \beta}(\beta ' + \a (b)) - 1 ) \leq e^{2 \beta} (\underbrace{|e^{- \beta} (\beta ' + \a (b))|}_{=: \Phi (\beta) (b)} - 1)\]

The second factor in the definition of $\Phi$ is a polynomial, so in the definition 2.2 from \cite{MN} we put $M= \R$, $\eta$ the Euclidean metric on $\R$, $k=1=d$, $\a =1$ and $u$ arbitrary. Then we define $P \in C^0(\R, \R{\rm Poly}_2^1): x \mapsto ((X,Y) \mapsto Y + \a (x))  $. Thus $\Phi$ is a flatzoomer (on $\R$). We define $\s_2:= \frac{1}{2} \ln (2 |\rho|/(n-2))$, then for $\beta > \s_2 $ we get $e^{2\beta} > e^{2 \s_2} = 2 |\rho|/(n-2)$. 

Moreover, we know by Eq. \ref{FZ1} that $\beta = \ln (\gamma/2) > ln (\eta) =:\s_1 $ implies future causal and b.a. completeness of $e^{2u} g$. We define $\s = \max \{ \s_1, \s_2\} $ pointwise. Then the main theorem Th 4.1 of \cite{MN} (with $M= \R$, $\e_i = 0$ for all $i \in \N$ and $w=f$) yields a smooth function $\beta > \s$ such that $\Phi (\beta) (b)< \frac{1}{2} \forall b \in \R$, which finally solves the problem, as then $\phi(\beta) (b) \leq 2 |\rho| (1/2-1)/(n-2) = - |\rho|/(n-2) < \rho/(n-2) $.

This finishes the proof of the spatially compact case.

\V For the general case, the only thing we need in order to transfer the proof above is properness of $t\vert_{J^+(S)}$, i.e., compactness of $t^{-1} ((-\infty; D)) \cap J^+(S)$ for every $D \in \R$. Indeed, we can construct such a 'future-proper' Cauchy temporal function by Theorem \ref{Lemma2}, which concludes the proof. \hfill \qed

\bigskip

\V {\bf Remark.} It is interesting that we could show by our method only that by appropriate conformal factors preserving the null energy condition we can reach causal completeness in one direction (future or past), whereas Bartnik's conjecture assumes two-sided causal completeness.

\newpage

{\bf Acknowledgements:} The author would like to thank Gregory Galloway, Eric Ling and two anonymous referees for helpful comments on previous versions of the article.

{\bf Data Availability Statement.} No experimental data has been produced linked to this article.

\end{document}